\documentclass[12pt]{article}
\usepackage{graphicx}
\usepackage{amsmath,amsthm,amssymb,enumerate}%, esint}
\usepackage{euscript,mathrsfs}
\usepackage{color}
\usepackage{dsfont}
\usepackage[left=2cm,right=2cm,top=3.5cm,bottom=3.5cm]{geometry}
\usepackage{color}
\usepackage[framemethod=tikz]{mdframed}
\usepackage{bm}
\allowdisplaybreaks

\usepackage{soul}

\catcode`\@=11 \@addtoreset{equation}{section}

\catcode`\@=12

\allowdisplaybreaks

\newtheorem{Theorem}{Theorem}[section]
\newtheorem{Proposition}[Theorem]{Proposition}
\newtheorem{Lemma}[Theorem]{Lemma}
\newtheorem{Corollary}[Theorem]{Corollary}

\theoremstyle{definition}
\newtheorem{Definition}[Theorem]{Definition}

\newtheorem{Remark}[Theorem]{Remark}

\newcommand{\bTheorem}[1]{
\begin{Theorem} \label{T#1} }
\newcommand{\eT}{\end{Theorem}}

\newcommand{\bProposition}[1]{
\begin{Proposition} \label{P#1}}
\newcommand{\eP}{\end{Proposition}}

\newcommand{\bLemma}[1]{
\begin{Lemma} \label{L#1} }
\newcommand{\eL}{\end{Lemma}}

\newcommand{\bCorollary}[1]{
\begin{Corollary} \label{C#1} }
\newcommand{\eC}{\end{Corollary}}

\newcommand{\bRemark}[1]{
\begin{Remark} \label{R#1} }
\newcommand{\eR}{\end{Remark}}

\newcommand{\bDefinition}[1]{
\begin{Definition} \label{D#1} }
\newcommand{\eD}{\end{Definition}}

\newcommand{\prst}{\mathbb{P}}

\newcommand{\jump}[1]{\left[ \left[ #1 \right] \right]}

\newcommand{\vrh}{\vr_h}

\newcommand{\tvm}{\widetilde{\vc{m}}}

\newcommand{\bfphi}{\boldsymbol{\varphi}}

\newcommand{\ds}{\,\mathrm{d}\sigma}

\newcommand{\Td}{\mathbb{T}^d}
\newcommand{\bFormula}[1]{
\begin{equation} \label{#1}}
\newcommand{\eF}{\end{equation}}

\newcommand{\grid}{\mathcal{T}}

\newcommand{\vuh}{\vu_h}

\newcommand{\intSh}[1] {\int_{\sigma} #1 \ds }

\newcommand{\Divh}{{\rm div}_h}
\newcommand{\Gradh}{\nabla_h}

\newcommand{\Ov}[1]{\overline{#1}}

\newcommand{\aleq}{\stackrel{<}{\sim}}

\newcommand{\avs}[1]{\left\{\!\!\left\{ #1\right\}\!\!\right\}}
\newcommand{\vvh}{\bm{v}_h}
\newcommand{\vQh}{\vc{Q}_h}
\newcommand{\pd}{\partial}

\newcommand{\gradd}{\nabla_{\mathcal D}}

\newcommand\Up{\mbox{\sl{Up}}}
\newcommand{\muh}{h^\varepsilon}

\newcommand{\norm}[1]{\left\lVert#1\right\rVert}

\newcommand{\vr}{\varrho}

\newcommand{\tvr}{\tilde \vr}
\newcommand{\tvu}{{\tilde \vu}}

\newcommand{\vu}{\vc{u}}
\newcommand{\vm}{\vc{m}}

\newcommand{\vn}{\vc{n}}

\newcommand{\vc}[1]{{\bm #1}}

\newcommand{\Div}{{\rm div}_x}
\newcommand{\Grad}{\nabla_x}

\newcommand{\dx}{\,{\rm d} {x}}

\newcommand{\dt}{\,{\rm d} t }

\newcommand{\intO}[1]{\int_{\Omega} #1 \ \D \mathbb{P}}

\newcommand{\intTd}[1]{\int_{\mathbb{T}^d} #1 \ \dx}

\newcommand{\vv}{\vc{v}}

\newcommand{\D}{{\rm d}}

\newcommand{\ep}{\varepsilon}

\newcommand{\expe}[1]{ \mathbb{E} \left[ #1 \right] }

\newcommand{\br}{ \nonumber \\ }

\def\softd{{\leavevmode\setbox1=\hbox{d}%
          \hbox to 1.05\wd1{d\kern-0.4ex{\char039}\hss}}}
\definecolor{Cgrey}{rgb}{0.85,0.85,0.85}
\definecolor{Cblue}{rgb}{0.50,0.85,0.85}
\definecolor{Cred}{rgb}{1,0,0}
\definecolor{fancy}{rgb}{0.10,0.85,0.10}

\newcommand\Cbox[2]{%
    \newbox\contentbox%
    \newbox\bkgdbox%
    \setbox\contentbox\hbox to \hsize{%
        \vtop{
            \kern\columnsep
            \hbox to \hsize{%
                \kern\columnsep%
                \advance\hsize by -2\columnsep%
                \setlength{\textwidth}{\hsize}%
                \vbox{
                    \parskip=\baselineskip
                    \parindent=0bp
                    #2
                }%
                \kern\columnsep%
            }%
            \kern\columnsep%
        }%
    }%
    \setbox\bkgdbox\vbox{
        \color{#1}
        \hrule width  \wd\contentbox %
               height \ht\contentbox %
               depth  \dp\contentbox
        \color{black}
    }%
    \wd\bkgdbox=0bp%
    \vbox{\hbox to \hsize{\box\bkgdbox\box\contentbox}}%
    \vskip\baselineskip%
}

\mdfdefinestyle{MyFrame}{%
	linecolor=black,
	outerlinewidth=1pt,
	roundcorner=5pt,
	innertopmargin=\baselineskip,
	innerbottommargin=\baselineskip,
	innerrightmargin=10pt,
	innerleftmargin=10pt,
	backgroundcolor=gray!20!white}

\date{}

\makeindex
\begin{document}

%%%%%%%%%%%%%%%%%%%%%%%%%%%%%%%%

\title{Convergence of a stochastic collocation finite volume method for the compressible Navier--Stokes system}

\author{Eduard Feireisl
	\thanks{The work of E.F. was partially supported by the
		Czech Sciences Foundation (GA\v CR), Grant Agreement
		18--05974S. The Institute of Mathematics of the Academy of Sciences of
		the Czech Republic is supported by RVO:67985840. \newline
		\hspace*{1em} $^\spadesuit$
		M.L. has been funded by the Deutsche Forschungsgemeinschaft (DFG, German Research Foundation) - Project number 233630050 - TRR 146 as well as by  TRR 165 Waves to Weather. She is grateful to the Gutenberg Research College
		and Mainz Institute of Multiscale Modelling for supporting her research.
	} \and  M\'aria Luk\'a\v{c}ov\'a-Medvi\softd ov\'a$^{\spadesuit}$
}

\date{\today}

\maketitle

\bigskip

\centerline{$^*$  Institute of Mathematics of the Academy of Sciences of the Czech Republic}

\centerline{\v Zitn\' a 25, CZ-115 67 Praha 1, Czech Republic}

\medskip

\centerline{$^\spadesuit$ Institute of Mathematics, Johannes Gutenberg--University Mainz}

\centerline{Staudingerweg 9, 55 128 Mainz, Germany}

\begin{abstract}
	
We propose a stochastic collocation method based on the piecewise constant interpolation on the probability space combined with a finite volume method to solve the compressible Navier--Stokes system at the nodal points. We show convergence of numerical solutions to a statistical solution of the Navier--Stokes system on condition that the numerical solutions are bounded in probability. The analysis uses the stochastic compactness method based on the Skorokhod/Jakubowski representation theorem and the criterion of convergence in probability due to Gy\" ongy and Krylov.

\end{abstract}


%{\bf MSC:}
\bigskip

%\tableofcontents

\section{Introduction}
\label{i}

Mathematical models arising in science and engineering inherit several sources of uncertainties, such as model parameters, initial and/or boundary conditions.
Prominent examples are the fluid flow models in meteorology, where the initial distribution of the {pressure (density)},
the velocity and the temperature, and even some rheological parameters, for instance the transport coefficients, may be viewed as random data, see e.g.~\cite{wiebe}.
Consequently, in order to predict {reliable} results, deterministic models are insufficient and more sophisticated methods are needed to analyse the influence of uncertainties on numerical solutions. In the recent years a wide variety of uncertainty quantification methods has been proposed and investigated. Although the standard Monte Carlo method  is often used in practical applications, it may become prohibitively expensive due to its slow convergence and large number of required samples. Alternatively, stochastic spectral methods, such as the stochastic Galerkin and stochastic collocation methods, can be applied in order to efficiently compute numerical solutions to systems with parametric uncertainty. Stochastic Galerkin method is based on a spectral element approximation in the probability space. It belongs to the class of intrusive methods, where the corresponding deterministic numerical scheme needs to be adjusted to calculate moment statistics of a solution. The stochastic collocation methods are non-intrusive and only require solving the underlying deterministic system at the certain collocation nodes together with application of a suitable interpolation method
in the probability space. Further details can be found in the monographs by  Le Ma\^{\i}tre and Knio~\cite{knio_book},
Pettersson et al.~\cite{pettersson_book}, Xiu~\cite{xiu_book}, Zhang and Karniadakis~\cite{karniadakis_book}.

Rigorous convergence analysis of these uncertainty quantification methods leans on uniqueness and
continuous dependence of solutions on random parameters (stochastic regularity of the solution). Relevant convergence results for the stochastic collocation methods can be found, for example, in Babu\v{s}ka et al.~\cite{babuska},  Nobile et al.~\cite{nobile}, Tang and Zhou \cite{tang} and the references therein. Convergence analysis of the stochastic Galerkin method was presented, e.g., by Babu\v{s}ka et al.~\cite{babuska2}, Bespalov et al.~\cite{bespalov}, Cohen et al.~\cite{cohen}, Ernst et al.~\cite{ernst}. For the convergence analysis of the Monte Carlo-type methods we refer to
Charrier et al.~\cite{Charrier}, Herrmann and Schwab~\cite{Herrmann},
Koley at al.~\cite{Koley}, Kuo et al.~\cite{Kuo}, Leonardi et al.~\cite{Leonardi},  Mishra and Schwab~\cite{Mishra_Schwab}.

Apparently, much less seems to be known for the nonlinear evolutionary equations arising in fluid flow modelling, where well posedness even in the deterministic setting represents a largely open problem.
Our aim is to show convergence of a stochastic collocation method for the random \emph{barotropic (isentropic) Navier--Stokes system}:

\begin{mdframed}[style=MyFrame]
	
	\begin{align}
		\partial_t \vr + \Div (\vr \vu) &= 0, \label{i1}\\
		\partial_t (\vr \vu) + \Div (\vr \vu \otimes \vu) + \Grad p(\vr) &= \Div \mathbb{S} (\Grad \vu),
		\label{i2} \\
		\mathbb{S}(\Grad \vu) &= \mu \left( \Grad \vu + \Grad^t \vu - \frac{2}{d} \Div \vu \mathbb{I} \right) +
		\eta \Div \vu \mathbb{I},\ \mu > 0,\ \eta \geq 0, \label{i3} \\
		p(\vr) &= a \vr^{\gamma},\ a > 0,\ \gamma > 1,
		\label{i4}
		\end{align}
	
	\end{mdframed}
describing the time evolution of the mass density $\vr = \vr(t,x)$ and the velocity field $\vu = \vu(t,x)$ of
a viscous, compressible barotropic fluid. For the sake of simplicity, we suppose the space periodic boundary conditions,
meaning the fluid domain may be identified with the flat torus,
\begin{equation} \label{i5}
\ t \in [0,T],\	x \in \Td, \ d=2,3.
	\end{equation}
Randomness is enforced through the initial data
\begin{mdframed}[style=MyFrame]
\begin{equation} \label{i6}
	\vr(0, \cdot) = \vr_0,\ \vm (0, \cdot) \equiv \vr \vu (0, \cdot) = \vm_0,
	\end{equation}

\end{mdframed}
and the viscosity coefficients $\mu$ and $\eta$ that
are random variables defined on a probability basis
\[
\Big( \Omega, \mathfrak{B}[\Omega], \prst \Big),
\]
where $\Omega$ is a compact metric space of events, $\mathfrak{B}[\Omega]$ the $\sigma-$field of Borel subsets of $\Omega$, and $\prst$ a complete Borel probability measure on $\Omega$.
To the best of our knowledge, this paper presents the first result available on convergence analysis of a numerical scheme approximating the random system \eqref{i1}--\eqref{i6}.

Opposite to the above mentioned results, the principal difficulty is due to the fact that the Navier--Stokes system \eqref{i1}--\eqref{i6} is not known to be solvable
in the class of smooth solutions on a possibly large time interval $(0,T)$
even for smooth initial data. The weak solutions exist globally in time for $\gamma > \frac{d}{2}$
(see \cite{EF70} and the pioneering work of Lions \cite{LI4}), however, their uniqueness
in terms of the initial data
remains an open problem. To approximate differential equations with low  regularity of exact solutions, several methods were proposed in the literature: The multi-element probabilistic collocation methods, see Foo and Karniadakis \cite{Foo}, Foo et al.~\cite{Foo2}, the
multi-resolution analysis methods using stochastic finite elements~Le Ma\^{\i}tre et al.~\cite{maitre} and multi-wavelet expansions Le Ma\^{\i}tre et al.~\cite{maitre2}.

Our stochastic collocation method for the random Navier-Stokes system is based on a piecewise constant interpolation of numerical solutions evaluated
at specific collocation nodes in the probability space and can be seen as a low order \emph{multi-element probabilistic collocation method}. Specifically, we consider a decomposition of the probabilistic basis $\Omega$,
\[
\Omega = \bigcup_{m=1}^{\nu(M)} \Omega^M_m,\ \Omega^M_m \in \mathfrak{B}[\Omega],\
\Omega_i^M \cap \Omega_j^M = \emptyset \ \mbox{for}\ i \ne j,
\]
and choose the nodal points
\[
\omega^M_m \in \Omega^M_m, \quad m=1,\dots, \nu(M).
\]
Given \emph{deterministic} initial data $(\vr^M_{0,m}, \vm^M_{0,m})$ evaluated at each nodal point,
\[
\vr^M_{0,m} (x) = \vr_0(x, \omega^M_m),\ \vm^M_{0,m}(x) = \vm_0(x, \omega_m^M),\ x \in \Td,
\]
along with the associated viscosity coefficients
\[
\mu_m^M = \mu(\omega_m^M),\ \eta_m^M = \eta(\omega_m^M),
\]
the exact solution $(\vr^M_m, \vu_m^M)_{m=1}^{\nu(M)}$ of the Navier--Stokes system will be approximated through a finite volume numerical scheme specified in Section \ref{FV} yielding a family of approximate solutions,
\[
(\vr_{h,m}^M, \vu_{h,m}^M), \ \mbox{where}\ h > 0 \ \mbox{denotes the mesh size.}
\]
For the sake of simplicity, we consider $h$ independent of $m$, however, $h = h(M)$,
\[
h(M) \to 0 \ \mbox{as}\ M \to \infty.
\]

This process gives rise to a sequence of discrete random variables $(\vr^M_h (t,x, \omega), $ $\vu^M_h (t,x, \omega))_{M=1}^\infty$
obtained using a suitable interpolation of the approximate solutions $(\vr_{h,m}^M, \vu_{h,m}^M),$ $m=1, \dots, \nu(M).$
Our goal is to show that for
\[
\max_{m \leq \nu(M)} \ {\rm diam}[ \Omega^M_m] \to 0,\ h(M) \to 0 \ \mbox{as}\ M \to \infty,
\]
the sequence of approximate solutions converges, specifically
\begin{equation} \label{i7}
\vr^M_{h(M)} \to \vr,\ \vu^M_{h(M)} \to \vu \ \ \mbox{as}\ M \to \infty\ \mbox{in probability},
\end{equation}
where $(\vr, \vu)$ is a statistical (random) solution to the Navier--Stokes system \eqref{i1}--\eqref{i6}, with the random initial data $(\vr_0, \vm_0)$ and the viscosity coefficients $\mu$ and $\eta$. For the time being, we leave unspecified the topology in which the convergence \eqref{i7} takes place.

To avoid the well--posedness problem in the class of weak solutions, we focus on
regular initial data, specifically,
\begin{align}
	\vr_0 &\in W^{3,2}(\Td),\ \inf_{x \in \Td} \vr_0(x) > 0, \label{i8}\\
	\vm_0 &\in W^{3,2}(\Td; R^d),\ \prst \mbox{-a.s.} \label{i9}
	\end{align}
Under these circumstances, the Navier--Stokes system admits a regular local in time solution, with a random life--span, see Matsumura and Nishida \cite{MANI}, Tani \cite{TAN}, among others. As shown by Sun, Wang, and Zhang \cite{SuWaZh1}, the eventual blow up of a smooth solution is conditioned by some concentration of the density. The leading idea of the present paper is therefore to exclude blow up of smooth solutions at least at a statistical level. Accordingly, our main working hypothesis is that the approximate solutions
$(\vr^{M}_{h(M)}, \vu^{M}_{h(M)})$ are \emph{bounded in probability}:
\begin{mdframed}[style=MyFrame]
	
\begin{align}
	\mbox{For any}\ \ep > 0, \ \mbox{there exists}\ N = N(\ep) \ \mbox{such that} \br
	\limsup_{M \to \infty} \prst \left[ \left( \| \vr^{M}_{h(M)} \|_{L^\infty((0,T) \times \Td)}
	 + \| \vu^{M}_{h(M)} \|_{L^\infty((0,T) \times \Td; R^d)} \right) > N
	\right] \leq \ep.
	\label{i10}
	\end{align}
\end{mdframed}

Hypothesis \eqref{i10} is quite mild as boundedness is required only for a statistical significant number of
approximate solutions. Indeed \eqref{i10} follows from a stronger hypothesis
\begin{equation} \label{i11}
	\expe{  \| \vr^{M}_{h(M)} \|_{L^\infty((0,T) \times \Td)}
		+ \| \vu^{M}_{h(M)} \|_{L^\infty((0,T) \times \Td; R^d)}  } \aleq 1,
\end{equation}
where $\mathbb{E}$ denotes the expected value with the respect to the measure $\prst.$
Note carefully that we do not impose any uniform bounds on the pointwise values (for fixed $\omega$) of the random variables $(\vr^{M}_{h(M)}, \vu^{M}_{h(M)})$.
We also remark that boundedness of a sequence of approximate solutions is anticipated frequently in analysis of many numerical methods.
Finally, boundedness in probability \eqref{i10} represents an
{\it a posteriori} condition verifiable directly in the course of a numerical simulation. In addition, we show that the satisfaction of \eqref{i10} for
a specific choice of collocation nodal points associated to a specific partition of the probability space implies convergence of the method for \emph{any} 	
choice of collocation nodal points and \emph{any} partition of $\Omega$.

Our approach is based on the stochastic compactness method, notably Skorokhod/Jakubowski representation
theorem developed in the context of weak (distributional) topologies in the monograph \cite{BrFeHobook}.
We proceed as follows:

\begin{itemize}
	\item Using a particular sample of random fields, the choice of which is motivated by \cite{BrFeHobook}, we use Skorokhod theorem to
	pass to a family of problems defined on a new probability space, where all quantities in question share the same law with the original ones. As an added benefit, the boundedness in probability stated in \eqref{i10} is transformed to a.s. boundedness on the new probability space.
	
	\item We perform the limit $M \to \infty$ on the new probability space. First we observe that
	the sequence of approximate numerical solutions converges, up to a subsequence,  to a dissipative measure--valued solution in the sense of \cite[Chapter 11]{FeLMMiSh} a.s. in the new probability space. Using the fact that the limit
	is bounded a.s. and the regularity criterion of Sun, Wang, and Zhang \cite{SuWaZh1} we conclude that the limit is a smooth solution of the Navier--Stokes system and the convergence is unconditional.
	
	\item We use Gy\" ongy--Krylov criterion to show convergence in probability of the approximate solutions on the original probability space.

	\end{itemize}

Finally, let us comment briefly on our choice of the stochastic collocation method. Here, we adopt the point of view that the real distribution of the data $(\vr_0, \vm_0)$ and the viscosity coefficients $\mu$, $\eta$ is {\it a priori} not known and the only piece of information available are their values at the nodes $\omega^M_m$. In contrast with the collocation methods based on a global polynomial approximation, we do not anticipate continuity of the random fields with respect to
the parameter $\omega$, which may be particularly relevant for the rheological parameters $\mu$ and $\eta$ that may attain only a specific \emph{finite} number of values. Instead we only assume that the piecewise constant interpolations
\begin{align}
(\vr_{0}^{M}, \vm_{0}^{M})(x, \omega) &= \sum_{m=1}^{\nu(M)} (\vr_0( x, \omega^M_m), \vm_0(x, \omega^M_m )) \mathds{1}_{\Omega^M_m}(\omega),\br (\mu^M, \eta^M)(\omega) &= \sum_{m=1}^{\nu(M)} (\mu(\omega^M_m), \eta(\omega^M_m)) \mathds{1}_{\Omega^M_m}(\omega),\ \nu(M) \to \infty \ \mbox{as}\ M \to \infty,
\label{i12}
\end{align}
approach, as $M \to \infty$, the initial data $(\vr_0, \vm_0)$
and the viscosity coefficients $(\mu, \eta)$
in probability for \mbox{any} choice of the
nodal points $\omega^M_m$ and the volumes $\Omega^M_m$. As we shall see below, this is equivalent to the assumption that the random variables in question are bounded and continuous (in $\omega$) with the exception of a
set of zero probability. Such a property is in fact equivalent to the Riemann integrability of the
random data in $\Omega$. In Section \ref{id}, we give a detailed proof  of this statement that may be of independent interest. Let us point out that
most of the real world application are based on
\emph{finite} dimensional probability space isomorphic to a cube. There is, however, a new rather
elegant theory of Riemann integration on general measured compact metric spaces
developed by Taylor \cite{Taylor}. We shall systematically refer to
Taylor's theory in the present paper. Note that such a choice of the probability space $\Omega$ goes beyond the standard Riemann/Darboux theory developed mostly on finite--dimensional spaces.

The ansatz \eqref{i12} corresponds to the piecewise constant interpolation obtained using (deterministic) numerical solutions,
\[
\vr^{M}_{h} (t,x, \omega) = \sum_{m=1}^{\nu(M)} \vr^{M}_{h,m}(t,x) \mathds{1}_{\Omega^M_m}(\omega),
\ \vu^{M}_{h}(t,x, \omega) = \sum_{m=1}^{\nu(M)} \vu^{M}_{h,m}(t,x) \mathds{1}_{\Omega^M_m}(\omega).
\]
Our goal is to show convergence of the resulting random
numerical solutions independently of the choice of the
collocation nodes under the sole assumption
\[
\max_{m = 1,\dots, \nu(M)} {\rm diam} [\Omega^M_m] \to 0 \ \mbox{as} \ M \to \infty.
\]
Note that our approximation method coincides with the \emph{nearest neighbour interpolation} in the case the partition of $\Omega$ is formed by a Voronoi tessellation.
We point out that our choice of this relatively simple interpolation method is motivated by the anticipated low regularity
of the data with respect to the random parameter.
The convergence of the sparse grids approximation  in the spirit of Smolyak~\cite{smolyak}, Xiu~\cite{xiu}, Nobile et al.~\cite{nobile} can be handled for the Navier--Stokes system in a similar manner.

The paper is organized as follows. In Section \ref{P}, we present the preliminary material including the exact formulation of the numerical scheme and state our main results. In Section \ref{id}, we show that
unconditional convergence of the piecewise constant data approximation is in fact equivalent to the Riemann integrability of the data. Section \ref{N} summarizes the properties of numerical solutions
obtained in the monograph \cite{FeLMMiSh}. In Section \ref{SC}, we apply the stochastic compactness method to transform the problem to a new probability space. The convergence proof is completed in Section \ref{G}
by means of the Gy\" ongy--Krylov criterion. The paper is concluded by a short discussion in Section~\ref{CRe}.

\section{Preliminaries, main results}
\label{P}
To state our main result in a rigorous way, we have to specify the probability basis as well as
the finite volume (FV) numerical method used to construct the approximate solutions.

\subsection{Probability basis}

Our choice of probability space $\Omega$ --a compact metric space-- is quite general and includes, in particular, the finite--dimensional case $\Omega \approx [0,1]^N$
with $\prst = \rho \mbox{d}y$, considered frequently in the literature, see e.g.  Babu\v{s}ka et al. \cite{babuska}.

\begin{Definition}[Partition] \label{DP1}
	
	A \emph{partition} of $\Omega$ is a (finite) family of Borel sets
	$( \Omega^M_m )_{m=1}^{\nu(M)}$,
	\[
	\Omega_i^M \cap \Omega_j^M = \emptyset \ \mbox{for}\ i \ne j,\ \Omega= \bigcup_{m=1}^{\nu(M)} \Omega^M_m.
	\]
	The \emph{diameter} of the partition $(\Omega_m^M)_{m=1}^{\nu(M)}$ is
\[
\mbox{diam}[(\Omega_m^M)_{m=1}^{\nu(M)}] = \max_{m=1, \dots, \nu(M)} \mbox{diam}[\Omega^M_m].
\]
The \emph{collocation nodes} associated to a partition $( \Omega^M_m )_{m=1}^{\nu(M)}$ are the points
$\omega^M_m$,
\[
\omega^M_m \in \Omega^M_m,\ m=1,\dots, \nu(M).
\]

	\end{Definition}

\subsection{Finite volume numerical scheme}
\label{FV}

We introduce a finite volume (FV) method to approximate the (deterministic) Navier-Stokes system \eqref{i1}-\eqref{i4}.
The physical domain $\Td$ is decomposed into finite volumes (cuboids for simplicity)
\[
 \Td  = \bigcup_{K \in \grid_h} K.
\]
The set of all faces $\sigma \in \partial K,$  $K \in \grid_h$ is denoted by $\Sigma.$
%while the set of faces on the %boundary is denoted by $\Sigma_{ext},$ and  the set of interior faces by  %$\Sigma_{int}=\Sigma\backslash \Sigma_{ext}$.
We suppose
$
|K| \approx h^d, \ |\sigma| \approx h^{d-1}
\ \mbox{ for any }\ K \in \grid_h, \mbox{ and } \sigma \in \Sigma,
$
where the parameter  $h\in(0,1)$  denotes the size of the mesh $\grid_h.$

The space of  functions constant on each element $K \in \grid_h$ is denoted ${Q}_h$. The associated
projection reads
\[
\Pi_h: L^1(\Td) \to {Q}_h,\ \Pi_h v = \sum_{K \in \grid_h} \mathds{1}_K \frac{1}{|K|} \int_K
v \dx.
\]
The average and jump of $v \in Q_h$ on a face $\sigma \in \Sigma$ are denoted
\[
\avs{v} = \frac{v^{\rm in} + v^{\rm out} }{2},\ \ \
\jump{ v }  = v^{\rm out} - v^{\rm in},
\]
where $v^{\rm out}, v^{\rm in}$ are respectively the outward, inward limits with respect to a given normal $\vn$ to $\sigma \in \Sigma.$

Moreover, the following discrete differential operators
for piecewise constant functions $r_h \in Q_h,$ $\vvh \in \vQh \equiv( Q_h)^d$ will be used:
\begin{equation*}
\begin{aligned}
&\Gradh r_h  = \sum_{K \in \grid_h}  (\Gradh r_h)_K 1_K,  \quad
(\Gradh r_h)_K = \sum_{\sigma\in \pd K} \frac{|\sigma|}{|K|} \avs{r_h} \vn,
\\
& \Delta_h r_h = \sum_{K\in \grid_h} (\Delta_h r_h)_K 1_K , \quad
 (\Delta_h r_h)_K = \sum_{\sigma \in \pd K}  \frac{|\sigma|}{|K|} \frac{\jump{r_h}}{h} ,
\\
& \gradd r_h =  \sum_{\sigma\in\pd K} \left(\gradd r_h \right)_{\sigma} 1_\sigma,  \quad \left(\gradd r_h\right) _{\sigma} =\frac{\jump{r_h} }{h} \vc{n}
\\
&\Divh \vvh  = \sum_{K \in \grid_h}  (\Divh  \vvh)_K 1_K, \quad
(\Divh \vvh)_K = \sum_{\sigma\in \pd K} \frac{|\sigma|}{|K|} \avs{\vvh} \cdot \vn.
\\
%&\Divup (r_h \vvh) = \sum_{K \in \grid_h}  1_K \Divup(r_h \vvh)_K, \quad
% \Divup(r_h \vvh)_K = \sumsfK \frac{|\sigma|}{|K|} F_h(r_h,\vvh).
 \end{aligned}
\end{equation*}

Convective terms will be approximated by a dissipative upwind numerical flux denoted by $F_h;$ specifically
\begin{eqnarray*}
\Up [r_h, \vv_h]  &=&
\avs{r_h} \ \avs{\vv_h} \cdot \vc{n} - \frac{1}{2} |\avs{\vv_h} \cdot \vc{n}| \jump{ r_h } \\
F_h(r_h,\vvh)
&=& {\Up}[r_h, \vvh] - \muh \jump{ r_h } = \avs{r_h} \ \avs{\vv_h} \cdot \vc{n}
- \left( \muh + \frac{1}{2} |\avs{\vv_h} \cdot \vc{n}| \right)\jump{ r_h },
\ -1 < \varepsilon .
\end{eqnarray*}
Analogously, we define the vector-valued numerical flux ${\bf F}_h({\bf r}_h,\vvh)$ componentwisely.

In order to discretize the time evolution in $[0,T]$ we introduce a time step $\Delta t > 0,$ $\Delta t \approx h,$ and denote
\[
t_k = k \Delta t,\ k=1,2,\dots, N_T.
\]
Furthermore, we set
\begin{align*}
v^k(x) = v(t^k,x)  \ \mbox{ for all } \  x\in \Td,\ t^k=k\,\Delta t \ \mbox{ for } k=0,1, \ldots, N_T.
\end{align*}
The time derivative $\frac{\partial {v} }{\partial t}$  is approximated by the backward Euler finite difference
\[
\frac{\partial {v} }{\partial t} \approx D_t {v}^k \equiv \frac{ {v}^k - {v}^{k-1} }{\Delta t}.
\]
Finally,
we introduce a piecewise constant interpolation in time of the discrete values $v^k$,
\begin{align}\label{NM_TD}
v_h(t,\cdot) = v_0 \mbox{ for } t<\Delta t,\ &
v_h(t,\cdot)=v^k \mbox{ for } t\in [k \Delta t,(k+1) \Delta t),\ k=1,2,\ldots,N_T.
\end{align}

\begin{Definition}[{\bf FV numerical scheme}] \label{DD1}
\phantom{mm}

\begin{itemize}
	\item Given the initial data $(\vr_{0}, \vm_{0}) \in L^1(\Td; R^{d+1})$, we set
\[
\vr_{0,h} = \Pi_h \vr_0,\ \vm_{0,h} = \Pi_h \vm_0, \ \vr_{0,h} \vu_{0,h} = \vm_{0,h}.
\]	

\item
\noindent A pair $(\vr_{h}, \vu_{h})$ of piecewise constant functions (in space and time) is a numerical approximation of the Navier-Stokes system \eqref{i1}-\eqref{i4}  if the following system of discrete equations holds:
\begin{subequations}\label{scheme}
\begin{align}
&\intTd{ D_t \vr_{h} \varphi_h } - \sum_{ \sigma \in \Sigma } \intSh{  F_h(\vr_{h},\vu_{h})
\jump{\varphi_h}   } = 0 \quad \mbox{for all } \varphi_h \in {Q}_h,\label{scheme_den}\\
&\intTd{ D_t  (\vr_{h} \vu_{h}) \cdot \bfphi_h } - \sum_{ \sigma \in \Sigma } \intSh{ {\bf F}_h(\vr_{h} \vu_{h},\vu_{h})
\cdot \jump{\bfphi_h}   }- \sum_{ \sigma \in \Sigma } \intSh{  \avs{p(\vr_{h})} \vc{n} \cdot \jump{ \bfphi_h }  } \nonumber \\
&= - \mu \frac{1}{h} \sum_{ \sigma \in \Sigma } \intSh{ \jump{\vu_{h}}  \cdot
\jump{\bfphi_h}  } - \lambda \intTd{ \Divh \vu_{h} \Divh \bfphi_h }
\quad \mbox{for all }
\bfphi_h \in \vQh . \label{scheme_mom}
\end{align}
\end{subequations}
where $\lambda =  \frac  1 d \mu +\eta$.

\end{itemize}
\end{Definition}

\begin{Definition}[{\bf Approximate statistical solution}] \label{stat_def}
	\phantom{mm}
	
\noindent	
Given a partition $\left(\Omega^M_m\right)_{m=1}^{\nu(M)}$  of $\Omega$ and a set of the collocation nodes $\omega^M_m \in \Omega^M_m$,
\emph{approximate statistical solution} of the (random) Navier--Stokes system is a pair of random variables,
		\[
		\vr^{M}_h(t,x, \omega) = \sum_{m=1}^{\nu(M)} \vr^M_{h,m} (t,x) \mathds{1}_{\Omega^M_m}(\omega),\
		\vu^M_h (t,x, \omega) = \sum_{m=1}^{\nu(M)} \vu^M_{h,m} (t,x) \mathds{1}_{\Omega^M_m}(\omega),
		\]
	 where
 $(\vr^{M}_{h,m}, \vu^M_{h,m})$ is a  solution of the FV method \eqref{scheme}, with the (deterministic) initial data $(\vr_0(\omega^M_m),\ \vm_0(\omega^M_m))$, and the viscosity coefficients $\mu(\omega^M_m)$, $\lambda (\Omega^M_m) = \frac{1}{d} \mu(\omega^M_m) + \eta(\omega^M_m)$.
\end{Definition}

\begin{Remark}
Statistical solutions given in Definition~\ref{stat_def} are sometimes called \emph{random solutions} in literature.
\end{Remark}

\subsection{Main result}

We start by introducing the class of admissible data. Let
$$
\mathcal{R}(\Omega, \prst) = \left\{f: \Omega \to R \ \Big|\ f \ \mbox{bounded},\   \prst\{ \omega \in \Omega\
\Big|\  f \mbox{ is not continuous at } \omega \} = 0 \right\}
$$
denote the class of Riemann integrable functions on $\Omega,$ cf. Taylor \cite{Taylor}. Moreover, we introduce the total energy as
\begin{equation}
E(\vr, \vm) \equiv \begin{cases}  \frac 1 2 \frac{|\vm|^2}{\vr} + P(\vr)    & \text{ if } \vr > 0,\\
                                     0 & \text{ if }  \vr = 0, \vm = 0 , \\
                                     \infty &  \text{ if } \vr = 0,  \vm \neq 0 \mbox{ or } \vr < 0,
\end{cases}
\end{equation}
where $P(\vr)$ is the pressure potential, $P'(\vr) \vr - P(\vr) = p(\vr).$ Specifically, we may consider $P(\vr) = \frac{a}{\gamma  - 1 } \vr^\gamma$ \, if \,  $p(\vr) = a \vr^\gamma.$
Note that $E$ is convex l.s.c. for $(\vr, \vm) \in R^{d+1}$. Alternatively, we also denote
$$
E(\vr, \vu) = \frac 1 2 \vr |\vu|^2 + P(\vr).
$$

\begin{Definition} [{\bf Admissible data}] \label{DM1}
	
	We say that the data $\vr_0$, $\vm_0$, $\mu$, $\eta$ are \emph{admissible} if:
	
	\begin{itemize}
		\item
	\begin{equation*}
	\vr_0(\cdot, \omega) \in W^{3,2}(\Td),\ \inf_{\Td} \vr_0 (\cdot, \omega) > 0,\ \vm_0(\cdot, \omega) \in W^{3,2}(\Td; R^3)\ \mbox{for}\ \omega \in \Omega;
	\end{equation*}
	\item
	\[
	\mu(\omega) \geq \underline{\mu} > 0, \ \eta(\omega) \geq 0 \ \mbox{for}\ \omega \in \Omega, \ \mbox{with  deterministic constant}\ \underline{\mu} > 0;
	\]	
	\item
	the functions
	\begin{align}
	\omega \in \Omega &\to \intTd{ \vr_0 (x, \omega) \varphi(x) },\ \varphi \in C^\infty(\Td), \br
	\omega \in \Omega &\to \intTd{ \vm_0 (x, \omega) \cdot \bfphi(x) },\ \bfphi \in C^\infty(\Td), \br
	\omega \in \Omega &\to \intTd{ E(\vr_0(x, \omega) , \vm_0(x, \omega)) }, \br
	\omega \in \Omega &\to \mu(\omega),\
	\omega \in \Omega \to \eta (\omega),
	\nonumber
	\end{align}
belong to $\mathcal{R}(\Omega, \prst).$
	
\end{itemize}	
	\end{Definition}

As we shall see in Section \ref{id}, admissibility of the data implies convergence in probability of their piecewise constant interpolation independent of the choice partition of the probability space and the position of the collocation nodes. Moreover, it can be shown that these two properties are in fact \emph{equivalent} provided the measure $\prst$ of any open ball in $\Omega$ and of its closure
coincides, see Beer \cite{Beer}. In view of this argument, admissibility in the sense of Definition \ref{DM1} is the weakest condition for a collocation method to be correctly defined.

Our main result reads as follows.

\begin{mdframed}[style=MyFrame]
	
		\begin{Theorem}[{\bf Convergence of stochastic collocation FV method}] \label{MT1}
			\phantom{mm}
			
			\noindent
		Let the initial data $(\vr_0, \vm_0)$ as well as the viscosity coefficients $\mu$, $\eta$ be admissible in the sense of Definition \ref{DM1}.
		
		Let $\left(\Omega^M_m \right)_{m=1}^{\nu(M)}$, $M=1,2, \dots$ be a sequence of partitions of $\Omega$ with a family of  nodal points
		$\omega^M_m \in \Omega^M_m$ such that
		\[
		\nu(M) \to \infty,\  {\rm diam}[(\Omega^M_m)_{m=1}^{\nu(M)}] \to 0 \ \mbox{as}\ M \to \infty.
		\]
		
		Suppose that for $h = h(M) \to 0$ as $M \to \infty,$ the associated sequence of
approximate statistical solutions introduced in Definition~\ref{stat_def} is bounded in probability:
		\begin{align}
			\mbox{For any}\ \ep > 0, \ \mbox{there exists}\ N = N(\ep) \ \mbox{such that} \br
			\limsup_{M \to \infty} \prst \left[ \left( \| \vr^{M}_{h(M)} \|_{L^\infty((0,T) \times \Td)}
			+ \| \vu^M_{h(M)} \|_{L^\infty((0,T) \times \Td; R^d)} \right) > N
			\right] \leq \ep. \label{hypothesis}
		\end{align}
		
		Then
		\begin{align}
			\vr^M_{h(M)} &\to \vr \ \ \mbox{as}\ M \to \infty \  \mbox{in}\ L^q((0,T) \times \Td),\ q < \infty, \ \mbox{in probability}, \br
			\vu^M_{h(M)} &\to \vu \ \ \mbox{as}\ M \to \infty \ \mbox{in}\ L^q((0,T) \times \Td; R^d),\ q < \infty, \ \mbox{in probability},
			\label{concl}
		\end{align}		
		where $(\vr, \vu)$ is the classical statistical solution of the Navier--Stokes system with the initial data
		$(\vr_0, \vm_0)$ and the viscosity coefficients $\mu$, $\eta$.		
	\end{Theorem}

	\end{mdframed}

\bigskip

\begin{Remark}[{\bf Convergence in expectation}] \label{MR1}
	
	As the initial energy $\intTd{ E(\vr_0, \vm_0) }$ is uniformly bounded in $\omega \in \Omega$, and the numerical solutions satisfy the energy inequality \eqref{N2bis} below, we have
	\[
	{\rm sup}_{t \in [0,T]} \| \vr^M_{h(M)} \|_{L^\gamma(\Td)}^\gamma + \int_0^T \intTd{ |\vu^{M}_{h(M)}|^2 } \aleq 1 \ \mbox{uniformly for}\ \omega \in \Omega.
	\]
Consequently, the convergence \eqref{concl} yields
\[
\expe{ \| \vr^{M}_{h(M)} - \vr \|_{L^\gamma((0,T) \times \Td)}^q } \to 0,\ \expe{ \| \vu^{M}_{h(M)} - \vu \|_{L^2((0,T) \times \Td; R^d)}^q } \to 0
\]
for any $1 \leq q < \infty$.	
	
	\end{Remark}

We point out that the above result is conditioned only by the satisfaction of \eqref{hypothesis}. In particular,
it is independent of the position of the collocation nodes and the choice of the ``elements''
$\Omega^M_m$. We also do not anticipate  continuity of the
data with respect to the random variable. As we have seen, hypothesis \eqref{hypothesis}
follows from a weaker stipulation \eqref{i11} that can be reformulated as

\begin{equation} \label{i11a}
\sum_{m=1}^{\nu(M)} \prst [\Omega^M_m] \left( \| \vr^M_{h(M),m} \|_{L^\infty((0,T) \times \Td)}
+ \| \vu^M_{h(M),m} \|_{L^\infty((0,T) \times \Td; R^d)} \right) \aleq 1 \ \mbox{uniformly for}\ M \to \infty,
\end{equation}
which can be explicitly controlled during numerical simulations.

Finally, we realize that the conclusion of Theorem \ref{MT1} implies that the Navier--Stokes
system admits a unique classical solution for the given set of random data $\prst-$a.s.
This yields the following corollary.

\begin{mdframed}[style=MyFrame]
	
	\begin{Corollary}[{\bf Unconditional convergence}] \label{CM1}
		
		Let the initial data $(\vr_0, \vm_0)$ as well as the viscosity coefficients $\mu$, $\eta$ be admissible in the sense of Definition \ref{DM1}. Suppose there
		there exists a sequence of partitions $\left(\Omega^M_m \right)_{m=1}^{\nu(M)}$ with a family of nodal points $\omega^M_m \in \Omega^M_m$ such that
		\[
		\nu(M) \to \infty,\ {\rm diam}[(\Omega^M_m)_{m=1}^{\nu(M)}] \to 0 \ \mbox{as}\ M \to \infty.
		\]
		In addition, suppose there is $h(M) \to \infty$ such that the associated sequence of
		approximate statistical solutions is bounded in probability in the sense of \eqref{hypothesis}.
		
		Then for any sequence of partitions $\left(\Omega^N_n \right)_{n=1}^{\nu(N)}$, $N=1,2,\dots$ satisfying
		\[
		\nu(N) \to \infty,\ {\rm diam}[(\Omega^N_n)_{n=1}^{\nu(N)}] \to 0 \ \mbox{as}\ N \to \infty,
		\]
		any family of nodal points $\omega^N_m \in \Omega^N_n$, and any $h(N) \to 0$ as $N \to \infty$,
		the associated family of approximate statistical solutions $(\vr^N_{h(N)}, \vu^N_{h(N)} )_{n=1}^{\nu(N)}$ converges, specifically,
		\begin{align}
			\vr^N_{h(N)} &\to \vr \ \mbox{as}\ N \to \infty \ \mbox{in}\ L^\gamma((0,T) \times \Td) \ \mbox{in probability}, \br
			\vu^N_{h(N)} &\to \vu \ \mbox{as}\ N \to \infty  \ \mbox{in}\ L^2((0,T) \times \Td; R^d) \ \mbox{in probability},
			\nonumber
		\end{align}		
		where $(\vr, \vu)$ is the classical statistical solution of the Navier--Stokes system with the initial data
		$(\vr_0, \vm_0)$ and the viscosity coefficients $\mu$, $\eta$.
		
		\end{Corollary}
	
	\end{mdframed}

The rest of the paper is devoted to the proof of Theorem~\ref{MT1} and Corollary \ref{CM1}.

\section{Data convergence}
\label{id}

We start with the following result on convergence of admissible data that can be of independent interest.

\begin{mdframed}[style=MyFrame]

\begin{Proposition}[{\bf Convergence of Riemann approximations}] \label{Pid1}
	Let $f \in \mathcal{R}(\Omega; \prst)$. Let $(\Omega^M_m)_{m=1}^{\nu(M)}$, $M=1,2,\dots$ be a sequence of partitions of $\Omega$ such that
	\[
	\nu(M) \to \infty,\ {\rm diam}[\left(\Omega^M_m \right)_{m=1}^{\nu(M)}] \to 0 \ \mbox{as}\ M \to \infty.
	\]
	Let
	\[
	f^M(\omega) = \sum_{m=1}^{\nu(M)} f(\omega^M_m) \mathds{1}_{\Omega^M_m}(\omega),\ \omega^M_m \in \Omega^M_m.
	\]
	Then
	\begin{equation} \label{idi1a}
	\expe{ \left| f^M - f \right|^q } \to 0 \ \mbox{as}\ M \to \infty
	\end{equation}
for any $1 \leq q < \infty$.	
	\end{Proposition}

\end{mdframed}

\begin{proof}
	
	We start by a result of Taylor \cite[Section 1, Section 3, Proposition 3.2 ]{Taylor}:
	\begin{align}
		f \in \mathcal{R}(\Omega; \prst) \ \Rightarrow \ f \ \mbox{is}\ \prst-\mbox{measurable, and} \br
	\intO{ f } = \lim_{M \to \infty} \sum_{m=1}^{\nu(M)} f(\omega^M_m) \prst[\Omega^M_m] = \lim_{M \to \infty} \intO{ f^M }.
	\label{id1}
	\end{align}
	
	Next, observe that
	\[
	(f^M)^2 = \left( \sum_{m=1}^{\nu(M)} f(\omega^M_m) \mathds{1}_{\Omega^M_m}(\omega) \right)^2 =
	\sum_{m=1}^{\nu(M)} f^2(\omega^M_m) \mathds{1}_{\Omega^M_m}(\omega).
		\]
	Seeing that
	\[
	f \in \mathcal{R}(\Omega; \prst) \ \Rightarrow f^2  \in \mathcal{R}(\Omega; \prst),
	\]
	we may apply \eqref{id1} to $f^2$ obtaining
	\[
	\intO{ (f^M)^2 } = \intO{ \sum_{m=1}^{\nu(M)} f^2(\omega^M_m) \mathds{1}_{\Omega^M_m}(\omega) } =
	\sum_{m=1}^{\nu(M)} f^2(\omega^M_m) \prst[\Omega^M_m] \to \intO{ f^2 }\ \mbox{ as } M\to \infty.
	\]
	
	Now, we examine the limit
	\[
	\lim_{M\to\infty}\intO{ f^M \phi },\ \phi \in {\rm Lip}(\Omega).
	\]
	Set
	\[
	\phi^M (\omega) = \sum_{m=1}^{\nu(M)} \phi (\omega^M_m) \mathds{1}_{\Omega^M_m}(\omega).
		\]
		Similarly to the above
	\[
	f^M \phi^M = \left( \sum_{m=1}^{\nu(M)} f(\omega^M_m) \mathds{1}_{\Omega^M_m}(\omega) \right)
	 \left( \sum_{m=1}^{\nu(M)} \phi(\omega^M_m) \mathds{1}_{\Omega^M_m}(\omega) \right) =
	\sum_{m=1}^{\nu(M)}  f(\omega^M_m)  \phi(\omega^M_m) \mathds{1}_{\Omega^M_m}(\omega).
	\]	
	As $f \phi \in \mathcal{R}(\Omega; \prst)$, we may infer, using \eqref{id1} once more,
	\begin{equation} \label{id2}
	\intO{ f^M \phi^M } \to \intO{ f \phi } \ \mbox{as}\ M \to
	\infty
	\end{equation}
	for any $\phi \in {\rm Lip}(\Omega)$. In addition, as $\phi$ is Lipschitz,
	\[
	\sup_{\Omega}|\phi^M - \phi| = \max_{m=1,\dots,\nu( M)} \sup_{\omega \in \Omega^M_m} |\phi (\omega^M_m) - \phi(\omega)| \leq L \max_{m=1,\dots, \nu(M)} {\rm diam}[\Omega^M_m] \to 0 \ \mbox{as}\ M \to \infty,
	\]
	where $L$ denotes the Lipschitz constant of $\phi$.
	As the sequence $f^M$ is uniformly bounded, the convergence \eqref{id2} yields
	\begin{equation} \label{id3}
		\intO{ f^M \phi } \to \intO{ f \phi } \ \mbox{as}\ M \to
		\infty \ \mbox{for any}\ \phi \in {\rm Lip}(\Omega).
	\end{equation}

Our ultimate goal is to extend \eqref{id3} to any $\phi \in L^1(\Omega; \prst)$. Since $(f^M)_{M=1}^\infty$ is uniformly bounded we need the set of Lipschitz functions to be dense in $L^1(\Omega; \prst)$. For measurable compact metric spaces
this was proved by Hanneke et al. \cite[Appendix A, Lemma A1]{Hanneke}.

Consequently, we have
\[
f^M \to f \ \mbox{weakly in}\ L^2(\Omega; \prst),\
\| f^M \|^2_{L^2(\Omega; \prst)} \to \| f \|^2_{L^2(\Omega; \prst)}
\]
yielding
\[
\expe{ | f^M - f |^2 } \to 0 \ \mbox{as}\ M \to \infty.
\]
As $f^M$ are uniformly bounded, the desired conclusion \eqref{idi1a} follows.

	\end{proof}

Since the data are admissible, Proposition \ref{Pid1} yields
\begin{align}
	\expe{ \left|  \intTd{\left(\sum_{m=1}^{\nu(M)} \vr_0 (\cdot, \omega^m_M ) \mathds{1}_{\Omega^m_M} - \vr_0 \right) \varphi } \right|^q } \to 0 \ \mbox{as}\ M \to \infty \ \mbox{for any}\ \varphi \in C^\infty(\Td), \br
	\expe{ \left|  \intTd{\left(\sum_{m=1}^{\nu(M)} \vm_0 (\cdot, \omega^m_M ) \mathds{1}_{\Omega^m_M} - \vm_0 \right) \bfphi }\right|^q } \to 0 \ \mbox{as}\ M \to \infty \ \mbox{for any}\ \bfphi \in C^\infty(\Td; R^3), \br
	\expe{ \left| \sum_{m=1}^{\nu(M)} \intTd{ E \Big(\vr_0 ( \cdot, \omega^m_M) ; \vm_0 ( \cdot, \omega^m_M) \Big) } \mathds{1}_{\Omega^m_M} -
		\intTd{ E ( \vr_0, \vm_0 )  }
		\right|^q } \to 0 \ \mbox{as}\ M \to \infty,
		\label{id4}
	\end{align}
and, similarly,
\begin{align}
\expe{ \left| \sum_{m=1}^{\nu(M)}  \mu(\omega^m_M)   \mathds{1}_{\Omega^m_M} -
		 \mu\right|^q } \to 0, \quad \expe{ \left| \sum_{m=1}^{\nu(M)}  \eta(\omega^m_M)   \mathds{1}_{\Omega^m_M} -
		 \eta\right|^q } \to 0
 \mbox{ as }\ M \to \infty, \label{id4a}
\end{align}
for any $1\leq q < \infty.$

Now we use the following result proved in \cite[Section~4.3]{MarEd}.
\begin{Lemma}\label{strong_conv}
Let
\begin{eqnarray*}
&&\vr_M \to \vr   \mbox{ weakly in }\ L^\gamma(\Td) \\
&&\vm_M \to \vm  \mbox{ weakly in }\ L^{\frac{2\gamma}{\gamma + 1}}(\Td; R^d)
\end{eqnarray*}
and
$$
\intTd{ E(\vr_M, \vm_M)}  \to \intTd{ E(\vr, \vm)} \mbox{ as } M \to \infty.
$$
Then
\begin{eqnarray*}
&&\vr_M \to \vr   \mbox{(strongly) in }\ L^\gamma(\Td) \\
&&\vm_M \to \vm  \mbox{ (strongly) in }\ L^{\frac{2\gamma}{\gamma + 1}}(\Td; R^d)  \mbox{ as } M \to \infty.
\end{eqnarray*}
\end{Lemma}

Obviously,
\begin{align}
 &\sum_{m=1}^{\nu(M)} \intTd{ E \Big(\vr_0 ( x, \omega^M_m) ; \vm_0 ( x, \omega^m_M) \Big) } \mathds{1}_{\Omega^M_m}(\omega) \br
 &= \intTd{ E \left( \sum_{m=1}^{\nu(M)} \vr_0 (x, \omega^M_m ) \mathds{1}_{\Omega^m_{\nu(M)}}(\omega); \sum_{m=1}^{\nu(M)} \vm_0 (x, \omega^M_m ) \mathds{1}_{\Omega^M_m} (\omega)   \right) },
 \nonumber
\end{align}
therefore it follows from \eqref{id4} and Lemma~\ref{strong_conv} that
\begin{align}
\sum_{m=1}^{\nu(M)} \vr_0 (\cdot, \omega^M_m ) \mathds{1}_{\Omega^M_m} &\to
\vr_0 \ \mbox{ in}\ L^\gamma(\Td)\ \prst-\mbox{a.s.}; \br
\sum_{m=1}^{\nu(M)} \vm_0 (\cdot, \omega^M_m ) \mathds{1}_{\Omega^M_m} &\to
\vm_0 \ \mbox{ in}\ L^{\frac{2\gamma}{\gamma + 1}}(\Td; R^3)\ \prst-\mbox{a.s.}
\nonumber
\end{align}
passing to a subsequence as the case may be.
As the energy
$$\intTd{E \left( \sum_{m=1}^{\nu(M)} \vr_0 (x, \omega^M_m ) \mathds{1}_{\Omega^M_m}(\omega); \sum_{m=1}^{\nu(M)} \vm_0 (x, \omega^M_m ) \mathds{1}_{\Omega^M_m} (\omega)   \right) }$$
is bounded,
we conclude
\begin{eqnarray}
\expe{ \left\| \sum_{m=1}^{\nu(M)} \vr_0 (\cdot, \omega^M_m ) \mathds{1}_{\Omega^M_m} - \vr_0 \right\|_{L^{\gamma}(\Td)} +
  \left\| \sum_{m=1}^{\nu(M)} \vm_0 (\cdot, \omega^M_m )\mathds{1}_{\Omega^M_m} -   \vm_0 \right\|_{L^{\frac{2\gamma}{\gamma+1}}(\Td; R^d)}}   \to 0  \label{id5}
\end{eqnarray}
as  $M \to \infty.$
Finally, as the limit is unique, the result is unconditional -- no need of subsequence.

\section{Properties of FV method}
\label{N}

 %Let
 %$(\vr_0, \vm_0)$ be given deterministic initial data. Let us denote the associated numerical solution
 %\[
 %\vrh, \vm_h, \ h > 0.
 %\]
%For the sake of simplicity, we consider the time step $\Delta t$ proportional to the parameter of the space discretization $h$. As the scheme is %time implicit, the numerical solutions may not uniquely determined by the initial data even et the discrete level. Fortunately, the choice is %restricted at each step to a finite number
%of initial data determined by the nodal points.

The crucial feature of FV method \eqref{scheme} is the fact that the approximate solutions $(\vrh, \vuh)$ give rise to a stable and consistent approximation of the Navier--Stokes system, see \cite{FeLMMiSh, FeiLukMizShe, FLM18_brenner}. Specifically, the following structure preserving properties hold:

\begin{itemize}
\item {\bf Conservation of discrete mass}\\
\begin{equation*}
\intTd{ \vr_{h} (\tau, \cdot) } = \intTd{ \vr_{0,h} } > 0,\,
\ \tau \geq 0.
\end{equation*}
\item {\bf Positivity  of the discrete density}\\
$$
\vr_{h} (\tau) >  0  \ \mbox{ for any } \tau > 0 \ \mbox{ provided }\
\vr_{0,h} > 0.
$$
\item {\bf Discrete total energy dissipation}\\
\begin{align*}
\intTd{ E( \vrh, \vuh ) (\tau, \cdot)  }
+  \int_0^T \mu   \norm{ \gradd  \vu_{h} }_{L^2(\Td)}^2  + \lambda \norm{\Divh \vu_{h} }_{L^2(\Td)}^2 \dt
&+\int_0^\tau \mathfrak{E}(\vrh, \vuh) \dt \br
&\leq
\intTd{ E( \vr_{0,h}, \vm_{0,h} )   },
\end{align*}
where the term $\mathfrak{E}(\vrh, \vuh) \geq 0$ represents numerical dissipation, see  \cite{FeLMMiSh} for its specific form.
\end{itemize}

\noindent As the energy is a convex function of $(\vr, \vm)$, Jensen's inequality yields
\begin{align}
\intTd{ E \left( \Pi_h [\vr_0], \Pi_h [\vm_0] \right) } &=
	\sum_{K \in \grid_h} \int_K E \left[ \frac{1}{|{K}|} \int_K \vr_0 (s) \D s, \frac{1}{|{K}|} \int_K \vm_0 (s) \D s \right] \mathds{1}_K \dx \br
	&\leq \sum_{K \in \grid_h} \int_K E \left( \vr_0, \vm_0 \right) \ \D s = \intTd{ E \left( \vr_0, \vm_0 \right) }.
	\nonumber
\end{align}
Thus the discrete energy inequality finally gives rise to
\begin{align}
	\intTd{ E( \vrh, \vuh ) (\tau, \cdot) } + \int_0^\tau \intTd{ \left[ \mu |\gradd \vuh |^2 +
		\lambda |{\rm div}_h \vuh |^2 \right] } \dt &+ \int_0^\tau \mathfrak{E}(\vrh, \vuh) \dt \br
	&\leq \intTd{ E( \vr_{0}, \vm_{0} ) },
	\label{N1}
\end{align}
where the right--hand side is independent of $h$.

\begin{itemize}
\item{\bf Consistent approximation}\\
As shown in \cite{FeLMMiSh, FeiLukMizShe} the numerical dissipation basically controls all consistency errors. Accordingly, the FV method \eqref{scheme} can be rewritten as follows
	\begin{align}
	\int_0^T& \intTd{ \left[ \vrh \partial_t \varphi + \vrh \vuh \cdot \Grad \varphi \right] } \dt =
	- \intTd{ \vr_{0,h} \varphi(0, \cdot) } + e_1[ h, \varphi ], \br
	\int_0^T &\intTd{ \left[ \vrh \vuh \cdot \partial_t \bfphi + \vrh \vuh \otimes \vu_h : \Grad \bfphi
		+ p(\vrh) \Div \bfphi \right] } \dt\br &= \int_0^T \intTd{ \mu \,\gradd \vuh : \Grad \bfphi } \dt
	+ \int_0^T \intTd{ \lambda \, {\rm div}_h \vuh {\rm div} \bfphi } \dt  \br  &-
	\intTd{ \vm_{0,h} \cdot \bfphi } + e_2[h, \bfphi ]
	\label{N2}
\end{align}
with the test functions $\varphi \in C^2([0,T] \times \Td), $ $\bfphi \in C^2([0,T] \times \Td; R^d),$ respectively.
The consistency errors are controlled in terms of the norm of the test functions and the initial energy.
Specifically, there is a function $B$ such that
\begin{align}
B : R^2 \to R \ \mbox{locally bounded}, \ o = o(h),\ o(h) \to 0 \ \mbox{as} \ h \to 0, \br
|e_{1}[h, \phi]| + |e_2[h, \bfphi] | \leq B \left( \| \varphi \|_{C^2([0,T] \times \Td)} + \| \bfphi \|_{C^2
([0,T] \times \Td; R^d)},\
\intTd{ E(\vr_0, \vm_0) } \right) o(h).
\label{N3}
\end{align}
\end{itemize}

\noindent Moreover, we write the initial data integral as
\[
\intTd{ \vr_{0,h} \phi } = \intTd{\vr_0 \phi } + \intTd{ (\Pi_h [\vr_0] - \vr_0) \phi },
\]
where
%\begin{align}
$$
\left| \intTd{ (\Pi_h [\vr_0] - \vr_0) \phi } \right| =
\left| \intTd{ (\Pi_h [\vr_0] - \vr_0) (\phi - \Pi_h [\phi ]) } \right| \br
\aleq h \| \vr_0 \|_{L^\gamma(\Td)} \| \phi \|_{C^1(\Td)}.
$$
%\end{align}
Similarly, we can control the initial data in the momentum equation. Accordingly, consistency equations
\eqref{N2} can be rewritten in the form
	\begin{align}
	\int_0^T& \intTd{ \left[ \vrh \partial_t \varphi + \vrh \vuh \cdot \Grad \varphi \right] } \dt =
	- \intTd{ \vr_{0} \varphi(0, \cdot) } + e_1[ h, \varphi ], \br
	\int_0^T &\intTd{ \left[ \vrh \vuh \cdot \partial_t \bfphi + \vrh \vuh \otimes \vu_h : \Grad \bfphi
		+ p(\vrh) \Div \bfphi \right] } \dt\br &= \int_0^T \intTd{ \mu \, \gradd \vuh : \Grad \bfphi } \dt
	+ \int_0^T \intTd{ \lambda \, {\rm div}_h \vuh {\rm div} \bfphi } \dt \br  & -
	\intTd{ \vm_{0} \cdot \bfphi(0, \cdot) } + e_2[h, \bfphi ]
	\label{N4}
\end{align}
with the consistency error estimate \eqref{N3} still valid.

\subsection{Dissipative solutions}
\label{d}

Any consistent approximation \eqref{N1}, \eqref{N4} generates a dissipative weak solution introduced in \cite{FeiLukMizShe}. We say that
$(\vr, \vu)$ is a  dissipative weak solutions of the Navier-Stokes system with the initial data $\vr_0$, $\vm_0$
if
\begin{align}
\vr &\in C_{\rm weak}([0,T]; L^\gamma(\Td)), \ \vr(0,\cdot) = \vr_0  \br	
\vu &\in L^2(0,T; W^{1,2}(\Td; R^d)),\br
\vr \vu &\in C_{\rm weak}([0,T]; L^{\frac{2 \gamma}{\gamma + 1}}(\Td; R^d)), \ \vr \vu(0, \cdot) = \vm_0
%\mathbb{S} &\in L^2(0,T; L^2(\Td; R^{d \times d}(\Td; R^d))
	\label{d1}
	\end{align}
and the following relations are satisfied in the sense of distributions:
	\begin{align}
		\partial_t \vr + \Div (\vr \vu) &= 0 ,
        \br
		\partial_t (\vr \vu) + \Div (\vr \vu \otimes \vu) + \Grad p(\vr) &= \Div \mathbb{S}(\Grad \vu) - \Div \mathfrak{R},
	  \br
		\intTd{ \left[ \frac{1}{2} \vr |\vu|^2 + P(\vr) \right] (\tau, \cdot)} + \Ov{d}
		\intTd{ {\rm tr}[\mathfrak{R}] (\tau, \cdot) } & +
		\int_0^\tau \intTd{ \left[ \mathbb{S}(\Grad \vu) : \Grad \vu \right] } \br
		&\leq \intTd{ \left[ \frac{1}{2} \frac{|\vm_0|^2}{\vr_0} + P(\vr_0) \right]},
		\label{d2}
	\end{align}
where $\Ov{d} > 0$ is a positive constant. The extra term
$\mathfrak{R} \in L^\infty(0,T;  \mathcal{M}^+(\Td; R_{\rm{sym}}^{d\times d}))$
called the Reynolds stress is a product of possible concentrations and/or oscillations in a generating sequence.
In the context of FV approximation, the relation $\mathfrak{R} = 0$ is essentially equivalent to strong convergence of $(\vr_h, \vr_h \vu_h)$ as $h \to 0$.

%\subsection{Conditional regularity}

The following conditional regularity shown in \cite[Corollary~6.4]{FeiLukMizShe} is essential:

\medskip
Suppose that the initial data belong to the class
\begin{equation}
\vr_0 \in W^{3,2}(\Td), \ \inf_{\Td} \vr_0 > 0, \ \ \vm_0 \in W^{3,2}(\Td; R^d), \label{init}
\end{equation}
and the dissipative weak solution satisfies
\[
\vr \in L^\infty((0,T) \times \Td), \ \ \vu \in L^\infty((0,T) \times \Td; R^d).
\]
Then  $(\vr, \vu)$ is a classical solution, and
\begin{equation} \label{d3}
	\mathfrak{R} \equiv  0.
\end{equation}

\subsection{Convergence of FV method}

Summarizing the above discussion we report the following result proved in \cite[Theorem~11.3]{FeiLukMizShe}.
	\begin{Proposition} [{\bf Convergence of FV method}] \label{PC1}
	
	Let $(\vrh, \vuh)_{h > 0}$ be a FV solution \eqref{scheme} with the initial data $(\vr_0, \vm_0),$
	\[
	{\rm ess} \inf_{\Td} \vr_0 > 0,\
	\intTd{ E(\vr_0, \vm_0) } < \infty.
	\]

Then the following holds:

\begin{enumerate}
	
	\item 	
	There exists $h_n \to 0$ such that
	\begin{align}
	\vr_{h_n} &\to \vr \ \mbox{weakly-(*) in}\ L^\infty(0,T; L^\gamma(\Td)), \ \br
	\vu_{h_n} &\to \vu \ \mbox{weakly in}\ L^2(0,T; L^r(\Td, R^d)),\ r = 6 \ \mbox{if}\ d=3,\ 1 \leq r < \infty\ \mbox{if}\ d = 2, \ \br
	\vm_{h_n} = \vr_{h_n} \vu_{h_n} &\to \vr \vu \ \mbox{weakly-(*) in} \ L^\infty(0,T; L^{\frac{2\gamma}{\gamma + 1}}(\Td; R^d)),
	\nonumber
	\end{align}
where $(\vr, \vu)$ is a dissipative weak solution of the Navier--Stokes system
specified in Section~\ref{d}.

\item

If, in addition,
\begin{itemize}
	\item the initial data $(\vr_0, \vu_0)$ belong to the regularity class \eqref{init},
	\item $\sup_{h > 0} \| \vrh \|_{L^\infty((0,T) \times \Td)} < \infty,$  \ \  $\sup_{h > 0} \| \vu_h \|_{L^\infty((0,T) \times \Td; R^d)} < \infty$,
	\end{itemize}
then the limit is the unique classical solution of the Navier--Stokes system and the convergence is unconditional, no need to subtract a subsequence. Moreover,
\begin{align}
 \vrh &\to \vr \ \mbox{in}\ L^q((0,T) \times \Td), \br	
 \vuh &\to \vu \ \mbox{in}\ L^{q}((0,T) \times \Td; R^d) \ \mbox{ for any}\ 1 \leq q < \infty.	\nonumber
	\end{align}

\item If the Navier--Stokes system admits a strong solution in the class
\begin{align}
\vr > 0, \ \vr &\in C([0,T] \times \Td), \ \Grad \vr,\ \vu \in C([0,T] \times \Td; R^d),  \br
\Grad \vu &\in C([0,T] \times \Td; R^{d \times d}),\ \partial_t \vu \in C([0,T] \times \Td; R^d),
\nonumber	
	\end{align}
then
\begin{align}
	\vrh &\to \vr \ \mbox{in}\ L^\gamma((0,T) \times \Td), \br	
	\vuh &\to \vu \ \mbox{in}\ L^2((0,T) \times \Td; R^d).	\nonumber
\end{align}

\end{enumerate}

	\end{Proposition}

\begin{Remark} \label{RC1}
	
	As a matter of fact, the conclusion of Proposition \ref{PC1} remains valid for any consistent approximation
	by piecewise constant functions satisfying \eqref{N4} together with the energy inequality \eqref{N1} perturbed by a small consistency error:
\begin{align}
	\intTd{ E( \vrh, \vuh ) (\tau, \cdot) } &+ \int_0^\tau \intTd{ \left[ \mu |\gradd \vuh |^2 +
		\lambda |{\rm div}_h \vuh |^2 \right] } \dt \br
	&\leq \intTd{ E( \vr_{0}, \vm_{0} ) } + e_3[h],\ e_3[h] \to 0 \ \mbox{as}\ h \to 0.
	\label{N2bis}
\end{align}	
	
	\end{Remark}

\section{Stochastic compactness method}
\label{SC}

Let $(\Omega^M_m)_{m=1}^{\nu(M)}$ be the partition considered in Theorem \ref{MT1}, with the corresponding sequence of
approximate statistical solutions $(\vr^M_{h(M)}, \vu^M_{h(M)})_{M=1}^{\infty}$,
with the viscosity coefficients
\[
\mu^M = \sum_{m=1}^{\nu(M)}  \mu( \omega^M_m )\mathds{1}_{\Omega^M_m},\ \ \eta^M = \sum_{m=1}^{\nu(M)} \eta(\omega^M_m)\mathds{1}_{\Omega^M_m},
\]
and the initial data
\[
\vr^M_0 = \sum_{m=1}^{\nu(M)} \vr_0(\omega^M_m) \mathds{1}_{\Omega^M_m},\
\vm^M_0 = \sum_{m=1}^{\nu(M)} \vm_0(\omega^M_m) \mathds{1}_{\Omega^M_m}.
\]
As we have shown in Section \ref{id}, formula \eqref{id5},
\begin{align}
	\expe{ \| \vr^M_0 - \vr_0 \|_{L^\gamma(\Td)} } &\to 0, \br
	\expe{ \| \vm^M_0 - \vm_0 \|_{L^{\frac{2\gamma}{\gamma + 1}}(\Td; R^d)} }	&\to 0, \br
	\expe{ \left| \intTd{ E(\vr^M_0, \vm^M_0 ) } - \intTd{ E(\vr_0, \vm_0) } \right| } &\to 0, \br
\expe{ \left| \mu^M -
		 \mu\right| } \to 0, \quad \expe{ \left| \eta^M  - \eta\right| } &\to 0\ \mbox{ as } M\to \infty.
\label{SC1}	
	\end{align}
In addition, as the initial data are admissible, cf.~\eqref{DM1}, in particular the initial energy is bounded, the discrete energy dissipation \eqref{N1} yields
\begin{align}
&\| \vr_{h(M)}^M \|_{L^\infty(0,T; L^\gamma(\Td))} + \| \vr_{h(M)}^M  \vu_{h(M)}^M\|_{L^\infty(0,T; L^{\frac{2\gamma}{\gamma + 1}}(\Td;R^d)} \aleq  1  \br
&\| \vu_{h(M)}^M \|_{L^2(0,T; L^r(\Td; R^d))} \aleq 1,\ r =6 \ \mbox{if}\ d = 3,\ 1 \leq r < \infty \ \mbox{if}\ d = 2,
\label{SC2}
\end{align}
for all  $M=1, \dots, $ and uniformly in  $\omega \in \Omega.$

\subsection{Application of Skorokhod's theorem}
\label{S}

Our goal is to convert hypothesis \eqref{hypothesis} stated in probability on $\Omega$ into a pointwise bound. To this end, we recall the celebrated Skorokhod theorem~\cite[Theorem~5.1, 5.2]{Bi}, \cite[Theorem~11.7.2]{dudley}.
\begin{Theorem}[{\bf Skorokhod theorem}] \label{TS1}
	
	Let $(\vc{U}^M)_{M=1}^\infty$ be a sequence of random variables ranging in a Polish space $X$. Suppose that
	their laws are tight in $X$, meaning for any $\ep > 0$, there exists a compact set $K(\ep) \subset X$
	such that
	\[
	\prst[ \vc{U}^M \in X \setminus K(\ep) ] \leq \ep \ \mbox{for all}\ M = 1,2,\dots.
	\]
	
	Then there is a subsequence
	$
	M_n \to \infty
	$
	and a sequence of random variables $(\widetilde{\vc{U}}^{M_n})_{n=1}^\infty$ defined on the standard probability space
	\[
	\Big(\widetilde{\Omega} = [0,1], \mathfrak{B}[0,1], \D y \Big)
	\]
	satisfying:
	\begin{itemize}
		\item
		\[
		\widetilde{\vc{U}}^{M_n} \approx_X \vc{U}^{M_n} \ \mbox{(they are equally distributed random variables)},
		\]
		\item
		\[
		\widetilde{\vc{U}}^{M_n} \to \widetilde{\vc{U}} \ \mbox{in}\  X \ \mbox{for every }  y \in [0,1].
		\]
	\end{itemize}
	
	\end{Theorem}
\begin{Remark}
The conclusion of  Theorem~\ref{TS1} are slightly stronger than the original version of Skorokhod and my be found in Jakubowski \cite{Jakub}.
\end{Remark}

Our goal is to apply the Skorokhod theorem to the sequence of random variables
\[
\vc{U}^M =
\left(\vr_0,\ \vr_0^M,\ \vm_0,\ \vm_0^M, \mu,\ \eta,\ \mu^M,\ \eta^M,\  \vr^M_{h(M)},\ \vu^M_{h(M)}, \Lambda^M \right),\ M =1,2,\dots
\]
where
\[
\Lambda^M =  \| (\vr^M_{h(M)}, \vu^M_{h(M)} ) \|_{L^\infty((0,T) \times \Td; R^4)}
\]
ranging in the Polish space
\[
X = (L^\gamma (\Td))^2 \times (L^{\frac{2 \gamma}{\gamma + 1}}(\Td; R^d))^2 \times R^4 \times W^{-k,2}((0,T) \times \Td) \times
W^{-k,2}((0,T) \times \Td; R^d) \times R,
\]
where $k \geq 4$ is an integer. In particular $W^{k,2}((0,T) \times \Td) \hookrightarrow \hookrightarrow C([0,T] \times \Td)$, $d=2,3$.
The negative Sobolev spaces could be replaced by more natural Lebesgue $L^p$ spaces endowed with the weak topology compatible with Jakubowski's
extension of Skorokhod's theorem, see \cite{Jakub}.

Now, tightness of the laws of $(\vr_0^M, \vm^M_0, \mu^M, \eta^M )_{M=1}^\infty$ in
$L^\gamma(\Td) \times L^{\frac{2 \gamma}{\gamma + 1}}(\Td; R^d) \times R^2$ follows from the convergence stated in \eqref{SC1}. Tightness of the laws of $(\vr^M_{h(M)},\ \vu^M_{h(M)})_{m=1}^M$ in
$W^{-k,2}((0,T) \times \Td) \times
W^{-k,2}((0,T) \times \Td; R^d)$
is a consequence of the bounds \eqref{SC2}. Finally, tightness of the laws of  $(\Lambda^M)_{M=1}^\infty$ in
$R$ follows from the hypothesis of boundedness in probability \eqref{hypothesis}. As a product of a final number of compacts is compact, we conclude that $(\vc{U}^M)_{M=1}^\infty$, $X$ complies with the hypotheses of
the Skorokhod theorem (Theorem \ref{TS1}).

The conclusion of Theorem \ref{TS1} can be stated in a simple form as the random variables in question are finitely distributed. Keeping the original labelling for a possible subsequence, we may infer that
for each $M$, there exists a partition of $[0,1]$ in $\nu(M)$ Lebesgue measurable sets $\widetilde{\Omega}_m^M \subset [0,1]$
such that the following holds:
\begin{itemize}
	\item
	\[
	\prst [\Omega^M_m] = | \widetilde{\Omega}^M_m |,\ m = 1,\dots, \nu(M);
	\]
	\item
	\[
	\tvr_0^M (x,y) = \sum_{m=1}^{\nu(M)} \vr_0 (x, \omega^M_m) \mathds{1}_{\widetilde{\Omega}^M_m}(y), \
	\tvm_0^M = \sum_{m=1}^{\nu(M)} \vm_0 (x, \omega^M_m) \mathds{1}_{\widetilde{\Omega}^M_m}(y),\
	x \in \Td,\ y \in [0,1];
	\]
	\item
	\begin{align}
	\widetilde{\vr}^{\nu(M)}_{h(M)}(t,x,y) &= \sum_{m=1}^M \vr^M_{h(m),m}(t,x) \mathds{1}_{\widetilde{\Omega}^M_m}(y),\
	\widetilde{\vu}^{\nu(M)}_{h(M)}(t,x,y) = \sum_{m=1}^M \vu^M_{h(m),m}(t,x) \mathds{1}_{\widetilde{\Omega}^M_m}(y), \br
	t &\in [0,T],\ x \in \Td,\ y \in [0,1].
	\label{SC3}
	\end{align}
	\end{itemize}
In addition, there are random variables
\[
\tvr_0 \in L^\infty(0,1; L^\gamma(\Td)), \ \tvm_0 \in L^\infty(0,1; L^{\frac{2 \gamma}{\gamma + 1}}(\Td; R^d)),\
\vr_0 \approx \tvr_0,\ \vm_0 \approx \tvm_0,
\]
\[
\tvr_0 \in W^{3,2}(\Td), \ \inf_{\Td} \tvr_0 > 0,\ \tvm_0 \in W^{3,2}(\Td; R^d) \ \mbox{for any}\ y \in [0,1]
\]
such that
\begin{align}
	\expe{ \left\| \tvr^M_0 - \tvr_0 \right\|_{L^\gamma(\Td)} } &\to 0 ,\br
		\expe{ \left\| \tvm^M_0 - \tvm_0 \right\|_{L^{\frac{2\gamma}{\gamma + 1}}(\Td; R^d)} } &\to 0, \br
		\expe{ \left| \intTd{ E(\tvr_0^M, \tvm_0^M) } - \intTd{ E(\tvr_0, \tvm_0) } \right| } &\to 0
		\ \mbox{as}\ M \to \infty.
	\label{SC4}
	\end{align}
Next,
\[
\widetilde{\mu}^{\nu(M)} = \sum_{m=1}^M \mu (\omega^M_m) \mathds{1}_{\widetilde{\Omega}^M_m}(y),\
\widetilde{\eta}^{\nu(M)} = \sum_{m=1}^M \eta (\omega^M_m) \mathds{1}_{\widetilde{\Omega}^M_m}(y),
\]
\begin{equation} \label{SC5}
\expe{ |\widetilde{\mu}^M - \widetilde{\mu} | } \to 0,\
\expe{ |\widetilde{\eta}^M - \widetilde{\eta} | } \to 0,\ \widetilde{\mu} \approx \mu,\
\widetilde{\eta} \approx \eta.
\end{equation}

Finally,
\[
\widetilde{\Lambda}^M =  \| (\tvr^M_{h(M)}, \tvu^M_{h(M)} ) \|_{L^\infty((0,T) \times \Td; R^4)}
\ \mbox{converges as}\ M \to \infty \ \mbox{for any}\ y,
\]
which yields the desired bound
\begin{equation} \label{SC6}
\| (\tvr^M_{h(M)}, \tvu^M_{h(M)} ) \|_{L^\infty((0,T) \times \Td; R^4)}	\
\mbox{bounded as}\ M \to \infty \ \mbox{for any}\ y \in [0,1].
	\end{equation}
Thus we may infer, combining \eqref{SC6} with the convergence claimed by Theorem \ref{TS1},
\begin{align}
	\tvr^M_0 &\to \tvr_0 \ \mbox{in}\ L^\gamma(\Td),\
	\tvm_0^M \to \tvm_0 \ \mbox{in}\ L^{\frac{2 \gamma }{\gamma + 1}}(\Td; R^d) \ \mbox{as}\ M \to \infty \  \mbox{for any}\ y \in [0,1], \br
	\widetilde{\mu}^M &\to \widetilde{\mu},\ 	\widetilde{\mu}^M \geq \underline{\mu} > 0,\ \widetilde{\eta}^M \to \widetilde{\eta} \  \mbox{as}\ M \to \infty \ \mbox{for any}\ y \in [0,1], \br
	\tvr^M_{h(M)} &\to \tvr \ \mbox{weakly-(*) in}\ L^\infty((0,T) \times \Td),\br \tvu^M_{h(M)} &\to \tvu \ \mbox{weakly-(*) in}\ L^\infty((0,T) \times \Td; R^d) \  \mbox{as}\ M \to \infty \ \mbox{for any}\ y \in [0,1].
\label{SC7}	
	\end{align}

\subsection{Limit system}

Our ultimate goal in this section is to identify the system of equations satisfied by the limit $(\tvr, \tvu)$
obtained in \eqref{SC7}. To this end, we fix $y \in [0,1]$, keeping in mind that each time we consider a subsequence, where the latter may depend on $y$. Fortunately, as we shall see below, the limit will be independent of
the choice of the subsequence therefore unconditional.

First, observe that \eqref{SC4} yields
\begin{equation} \label{SC8}
\intTd{ E(\tvr^M_0, \tvm^M_0 ) } \to \intTd{ E(\tvr_0, \tvm_0 ) } \ \mbox{as}\ M \to \infty\  \mbox{a.s. in}\ [0,1]
\end{equation}
at least for a suitable subsequence.

In accordance with \eqref{SC3}, the pair $(\tvr^M_{h(M)}, \tvu^M_{h(M)})$ is a solution of FV scheme \eqref{scheme} with the initial data $(\tvr^M_0, \tvm^M_0)$ for any $y \in [0,1]$. Specifically, the
consistency formulation \eqref{N1}, \eqref{N4},
\begin{align}
	\intTd{ E( \tvr^M_{h(M)}, \tvu^M_{h(M)} ) (\tau, \cdot) } &+ \int_0^\tau \intTd{ \left[ \widetilde{\mu}^M |\gradd \tvu^M_{h(M)} |^2 +
		\widetilde{\lambda}^M |{\rm div}_h \tvu^M_{h(M)} |^2 \right] } \dt  \br
	&\leq \intTd{ E( \tvr^M_{0}, \tvm^M_{0} ) },
	\label{SC9}
\end{align}
	\begin{align}
	\int_0^T& \intTd{ \left[ \tvr^M_{h(M)} \partial_t \varphi + \tvr^M_{h(M)}  \tvu^M_{h(M)} \cdot \Grad \varphi \right] } \dt =
	- \intTd{ \tvr^M_{0} \varphi(0, \cdot) } + e_1[ h(M) , \varphi ], \br
	\int_0^T &\intTd{ \left[ \tvr^M_{h(M)}  \tvu^M_{h(M)} \cdot \partial_t \bfphi + \tvr^M_{h(M)}  \tvu^M_{h(M)} \otimes  \tvu^M_{h(M)} : \Grad \bfphi
		+ p(\tvr^M_{h(M)}) \Div \bfphi \right] } \dt\br &= \int_0^T \intTd{ \widetilde{\mu}^M \, \gradd  \tvu^M_{h(M)} : \Grad \bfphi } \dt
	+ \int_0^T \intTd{ \widetilde{\lambda}^M \, {\rm div}_h  \tvu^M_{h(M)} {\rm div} \bfphi } \dt \br  & -
	\intTd{ \tvm^M_{0} \cdot \bfphi(0, \cdot) } + e_2[h(M), \bfphi ],
	\label{SC10}
\end{align}
holds for any $y \in [0,1]$.

Now, we use the fact that the viscosity coefficients $\widetilde{\mu}^M$ are bounded below by a positive deterministic constant $\underline{\mu}$ to deduce from the energy inequality \eqref{SC9} a uniform bound
\begin{equation} \label{SC11}
	\int_0^T \intTd{ \left[ |\gradd \tvu^M_{h(M)} |^2 +
		 |{\rm div}_h \tvu^M_{h(M)} |^2 \right] } \dt \aleq 1 \ \mbox{uniformly for}\ y \in [0,1].
	 \end{equation}
This allows us to rewrite the consistency formulation \eqref{SC9}, \eqref{SC10} in the form:
\begin{align}
	\intTd{ E( \tvr^M_{h(M)}, \tvu^M_{h(M)} ) (\tau, \cdot) } &+ \int_0^\tau \intTd{ \left[ \widetilde{\mu} |\gradd \tvu^M_{h(M)} |^2 +
		\widetilde{\lambda} |{\rm div}_h \tvu^M_{h(M)} |^2 \right] } \dt  \br
	&\leq \intTd{ E(\tvr_0, \tvm_0) } \br
	&+ \int_0^\tau \intTd{ \left[ (\widetilde{\mu} - \widetilde{\mu}^M) |\gradd \tvu^M_{h(M)} |^2 +
		(\widetilde{\lambda} - \widetilde{\lambda}^M) |{\rm div}_h \tvu^M_{h(M)} |^2 \right] } \dt \br &+
	\intTd{ E( \tvr^M_{0}, \tvm^M_{0} ) - E(\tvr_0, \tvm_0) }, \br
	\label{SC12}
\end{align}
where, by virtue of the bound \eqref{SC11} and the convergence stated in \eqref{SC7}, \eqref{SC8},
\begin{align}
e_3[h(M)] &\equiv \int_0^T \intTd{ \left[ (\widetilde{\mu} - \widetilde{\mu}^M) |\gradd \tvu^M_{h(M)} |^2 +
	(\widetilde{\lambda} - \widetilde{\lambda}^M) |{\rm div}_h \tvu^M_{h(M)} |^2 \right] } \dt \br
&+ \intTd{ E( \tvr^M_{0}, \tvm^M_{0} ) - E(\tvr_0, \tvm_0) } \to 0 \ \mbox{as}\ M \to \infty
\ \mbox{a.s. in} \ [0,1].
\nonumber
\end{align}

Using similar arguments, we can deduce from \eqref{SC10}:
	\begin{align}
	\int_0^T& \intTd{ \left[ \tvr^M_{h(M)} \partial_t \varphi + \tvr^M_{h(M)}  \tvu^M_{h(M)} \cdot \Grad \varphi \right] } \dt =
	- \intTd{ \tvr_{0} \varphi(0, \cdot) } + e_1[ h(M) , \varphi ], \br
	\int_0^T &\intTd{ \left[ \tvr^M_{h(M)}  \tvu^M_{h(M)} \cdot \partial_t \bfphi + \tvr^M_{h(M)}  \tvu^M_{h(M)} \otimes  \tvu^M_{h(M)} : \Grad \bfphi
		+ p(\tvr^M_{h(M)}) \Div \bfphi \right] } \dt\br &= \int_0^T \intTd{ \widetilde{\mu} \, \gradd  \tvu^M_{h(M)} : \Grad \bfphi } \dt
	+ \int_0^T \intTd{ \widetilde{\lambda} \, {\rm div}_h  \tvu^M_{h(M)} {\rm div} \bfphi } \dt \br  & -
	\intTd{ \tvm_{0} \cdot \bfphi(0, \cdot) } + e_2[h(M), \bfphi ] \quad  \mbox{a.s. in}\ [0,1],
	\label{SC13}
\end{align}
where the errors resulting from replacing $(\tvr^M_0, \tvm^M_0)$ by $(\tvr_0, \tvm_0)$, and
$(\widetilde{\mu}^M, \widetilde{\eta}^M)$ by $(\widetilde{\mu}, \widetilde{\eta})$ can be incorporated
in $e_1$, $e_2$. Thus, we may infer that $(\tvr^M_{h(M)} , \tvu^M_{h(M)})_{M=1}^\infty$ is a consistent
approximation of the Navier--Stokes system with the viscosity coefficients $\widetilde{\mu}$, $\widetilde{\lambda}$ and the initial data
$(\tvr_0, \tvm_0)$ in the sense specified in \eqref{N4}, \eqref{N2bis} for a.a. $ y \in [0,1]$.

In view of \eqref{SC6}, the approximate solutions are uniformly bounded; whence we may use the convergence result
stated in Proposition \ref{PC1} to conclude that
\begin{eqnarray*}
&&\tvr^M_{h(M)} \to \tvr \quad \mbox{ in } L^q((0,T) \times \Td) \\
&& \tvu^M_{h(M)} \to \tvu \quad \mbox{ in } L^q((0,T) \times \Td; R^d) \quad \mbox{ a.s. in } [0,1],\ 1 \leq q < \infty,
\end{eqnarray*}
where $(\tvr, \tvu)$ is the unique classical solution of the Navier-Stokes system with the initial data
$(\tvr_0, \tvm_0)$ and viscosity coefficients $\widetilde{\mu}, \widetilde{\lambda}.$
As the limit is unique, the convergence is unconditional.

\section{Application of Gy\" ongy -- Krylov theorem}
\label{G}

The ultimate goal is to show the convergence of the numerical method in terms of the original probability space. To this end we apply the following result, see Gy\" ongy, Krylov~\cite[Lemma~1.1]{Gkrylov}.

\begin{Theorem}[{\bf Gy\" ongy--Krylov theorem}]
	
	Let $X$ be a Polish space and $(\vc{U}^M)_{M \geq 1}$ a sequence of $X-$valued random variables.
	
	Then $(\vc{U}^M )_{M = 1}^\infty$ converges in probability if and only if for any sequence of joint laws of
	\[
	(\vc{U}^{M_k}, \vc{U}^{N_k})_{k = 1}^\infty
	\]
	there exists further subsequence that converge weakly to a probability measure $\mu$ on $X \times X$ such that
	\[
	\mu \left[ (x,y) \in X \times X,\ x = y \right] = 1.
	\]
	\end{Theorem}

Under the hypotheses of the previous section, we consider the sequence
\begin{align}
\left(\vr_0,\ \vr_0^{M_k},\ \vm_0,\ \vm_0^{M_k}, \mu,\ \eta,\ \mu^{M_k},\ \eta^{M_k},\  \vr^{M_k}_{h(M_k)},\ \vu^{M_k}_{h(M_k)},
\right.
\ \Lambda^{M_k},
\br
\left.
\vr_0^{N_k}, \ \vm_0^{N_k},\ \mu^{N_k},\ \eta^{N_k},\ \vr^{N_k}_{h(N_k)}, \ \vu^{N_k}_{h(N_k)},\  \Lambda^{N_k}\right)_{k=1}^\infty.
\nonumber
\end{align}
Similarly to the preceding section, we obtain the Skorokhod representation
\[
\tvr^{M_k}_{h(M_k)}, \ \tvr^{N_k}_{h(N_k)},\ \tvm^{M_k}_{h(M_k)},\ \tvm^{N_k}_{h(N_k)}
\]
satisfying
\begin{align}
(\tvr^{M_k}_{h(M_k)}, \tvr^{N_k}_{h(N_k)}) &\to (\tvr, \tvr) \ \mbox{in}\ L^q((0,T) \times \Td) \ \mbox{ a.s. in } [0,1], \br
(\tvu^{M_k}_{h(M_k)}, \tvu^{N_k}_{h(N_k)}) &\to (\tvu, \tvu) \ \mbox{in}\ L^q((0,T) \times \Td; R^3) \ \mbox{ a.s. in }[0,1],\ 1\leq q < \infty,
\end{align}
where $(\tvr, \tvu)$ is the unique classical solution of the Navier-Stokes system with the initial data
$(\tvr_0, \tvm_0)$ and viscosity coefficients $\widetilde{\mu}, \widetilde{\lambda}.$

Thus applying the Gy\" ongy--Krylov theorem to the sequence
\[
(\vr^M_{h(M)}, \vu^M_{h(M)})_{M=1}^\infty \ \mbox{ranging in the Polish space}\ L^q(\Td; R^4),
\]
we may infer that
\begin{equation} \label{G1}
	\vr^M_{h(M)} \to \vr, \ \vu^M_{h(M)} \to \vu \ \mbox{ in } L^q(\Td; R^4), \ 1\leq q < \infty, \mbox{ in probability,}
	\end{equation}
where $(\vr, \vu)$ is the unique solution of the Navier--Stokes system with the initial data
$(\vr_0, \vm_0)$ and the viscosity coefficients ${\mu}, {\eta}.$

We have shown Theorem~\ref{MT1}.

\subsection{Unconditional convergence}

It remains to show Corollary \ref{CM1}. As a byproduct of the conclusion of Theorem \ref{MT1}, we already know that the Navier--Stokes system admits a classical solution in $([0,T] \times \Td)$
for $\prst-$ a.a. data. Consequently, Corollary \ref{CM1} follows directly from Proposition \ref{PC1}, part 3. Indeed the convergence is now unconditional and there is no need of the stochastic compactness method.

\section{Concluding remarks}
\label{CRe}

In this paper we have studied convergence of a stochastic collocation FV method for the random compressible Navier--Stokes system. We have clarified the necessary conditions for the stochastic collocation to be
meaningful, meaning arbitrary choice of the collocation points and the ``elements'' in the probability space gives rise to the same asymptotic limit. We also eliminated the ambiguity in the limit of approximate numerical solutions by requiring that the approximate solutions are bounded in probability. The result can be seen
as the first rigorous proof of convergence of a statistical method for the compressible Navier--Stokes system.

We have made several constitutive restrictions that can be easily removed. In particular, the specific isentropic
form of the pressure is not really necessary, any EOS preserving monotonicity of $p$ as a functions of $\vr$ and the asymptotic behaviour $p(\vr) \approx \vr^\gamma$ for $\vr \to \infty$ can be handled by the same method.

The Navier--Stokes system can be augmented by a driving force that can be random similarly to the viscosity coefficients $\mu$ and $\eta$.

The results actually holds for any consistent approximation of the Navier--Stokes system not necessarily generated by a numerical method.

Finally, let us comment shortly on \emph{boundedness in probability} postulated in hypothesis \eqref{hypothesis}. The situation is particularly simple if we assume some
pointwise bounds though not on a very large set. As we have observed in Corollary \ref{CM1}, it is enough to find \emph{one} particular
sequence of partitions $(\Omega^M_m )_{m=1}^{\nu(M)}$, $M=1,2,\dots$, with the associated collocation nodes $\omega^M_m$ such that \eqref{hypothesis} holds to guarantee convergence of the method
for \emph{any} sequence of partitions. The following result (Taylor \cite[Proposition 5.1]{Taylor}) comes handy.

\begin{Lemma} \label{LC1}
	For any $\ep > 0$, there exists a partition $(\Omega^M_m)_{m=1}^{\nu(M)}$ such that
	\[
	{\rm diam}[ \Omega^M_m] < \ep,\ \prst[ \partial \Omega^M_m] = 0 \ \mbox{for all}\ m=1, \dots, \nu(M).
	\]
	\end{Lemma}

Now, given $K > 0$ consider the set $R(K) \subset \Omega$ of ``regular events'',
\begin{align}
R(K) = &\left\{ \omega \in \Omega \ \Big| \  \ \mbox{there exists}\ \Ov{h} > 0
\ \mbox{and a family of FV solutions}\  (\vr_{h}, \vu_{h} )_{0 < h < \Ov{h} } \right. \br
&\quad \mbox{with the initial data}\ (\vr_0, \vm_0)(\omega) \ \mbox{and the viscosity coefficients}
 (\mu, \eta)(\omega) \br
&\quad \mbox{such that}
\|(\vr_{h}, \vu_{h} )\|_{L^\infty(\Td; R^{d + 1})} \leq K \ \mbox{for all}\ 0 < h < \Ov{h} \Big\}.
\nonumber
\end{align}
Consider the following hypothesis:
\begin{equation} \label{C1}
	\ \mbox{There exists}\ K > 0 \ \mbox{such that}\ R(K) \ \mbox{is dense in}\ \Omega.
	\end{equation}
Note that a dense set may not be ``big'' in the measure sense, in particular, we may have
\[
\prst[ R(K) ] = 0.
\]
We claim that \eqref{C1} yields boundedness in probability for a suitable choice of partitions. Indeed, for any $\ep = \ep(M) \to 0$, Lemma \ref{LC1} yields a partition
$(\Omega^M_m)_{m=1}^{\nu(M)}$ such that
\[
\sup_{m = 1, \dots, \nu(M) } {\rm diam}[ \Omega^M_m] < \ep (M),\ \prst[ \partial \Omega^M_m] = 0 \ \mbox{for all}\ m=1, \dots, \nu(M).
\]
In particular,
\[
\prst[ \Omega^M_m ] = \prst[ {\rm int} [\Omega^M_m] ].
\]
Now, we choose the collocation nodes:
\[
\omega^M_m \in R(K) \ \mbox{if} \ {\rm int}[\Omega^M_m] \ne \emptyset,\ \omega^M_m \ \mbox{arbitrary otherwise.}
\]
Moreover, as the partition is finite, we may fix $h = h(M)$ so that
\[
\|(\vr^M_{h(M)}, \vu^M_{h(M)} )\|_{L^\infty(\Td; R^{d + 1})} \leq K \ \ \prst-\mbox{a.s.}
\]
As this construction can be repeated for any $M = \ep (M) \to 0$, the desired conclusion follows.

\def\cprime{$'$} \def\ocirc#1{\ifmmode\setbox0=\hbox{$#1$}\dimen0=\ht0
	\advance\dimen0 by1pt\rlap{\hbox to\wd0{\hss\raise\dimen0
			\hbox{\hskip.2em$\scriptscriptstyle\circ$}\hss}}#1\else {\accent"17 #1}\fi}

\end{document}